\documentclass[english]{amsart}
\usepackage[T1]{fontenc}
\usepackage[latin1]{inputenc}
\usepackage{babel}
\usepackage{amssymb}

\makeatletter

 \theoremstyle{plain}    
 \newtheorem{thm}{Theorem}[section]
 \numberwithin{equation}{section} 
 \numberwithin{figure}{section} 
 \theoremstyle{plain}    
 \newtheorem{cor}[thm]{Corollary} 
 \theoremstyle{plain}    
 \newtheorem{lem}[thm]{Lemma} 
 \theoremstyle{plain}    
 \newtheorem{prop}[thm]{Proposition} 
 \theoremstyle{definition}
 \newtheorem{defn}[thm]{Definition}
 \theoremstyle{definition}
  \newtheorem{example}[thm]{Example}
 \theoremstyle{remark}
 \newtheorem{rem}[thm]{Remark}
 \theoremstyle{remark}    
 \newtheorem{notation}[thm]{Notation} 

\usepackage[vcentermath,enableskew]{youngtab}
\usepackage[dvips]{graphicx}
\usepackage[dvips]{color}

\makeatother
\begin{document}

\title{Quantum cohomology of Grassmannians modulo symmetries}

\author{Harald Hengelbrock}



\curraddr{Fakultät für Mathematik, Ruhr-Universität Bochum, 44803 Bochum}

\email{harald.hengelbrock@ruhr-uni-bochum.de}

\date{November 4, 2002}

\begin{abstract}
The quantum cohomology of Grassmannians exhibits two symmetries related
to the quantum product, namely a \( \Bbb {Z}/n \) action and an involution
related to complex conjugation. We construct a new ring by dividing
out these symmetries in an ideal theoretic way and analyze its structure,
which is shown to control the sum of all coefficients appearing in
the product of cohomology classes. We derive a combinatorial formula
for the sum of all Littlewood-Richardson coefficients appearing in
the expansion of a product of two Schur polynomials.
\end{abstract}
\maketitle

\section{Introduction}

In the quantum Schubert calculus of the Grassmannians there appear
two kinds of important symmetries, one a cyclic symmetry as described
in \cite{AW}, another is an involution related to complex conjugation,
see \cite{Hen,Pos,Rie}. In the representation of Schubert classes
as complex functions on a finite set of points in complex space, these
symmetries act by changing the {}``complex direction'' of the function
values. In order to understand the invariant aspects of quantum Schubert
calculus it should therefore be rewarding to understand the absolute
value of the Schubert class functions. It turns out that the value
of Schubert classes on certain real points already determines an interesting
numerical aspect of this structure. This leads to the definition of
a quotient ring, which will be the coordinate ring of these points.
Next, the quantum Schubert calculus modulo symmetries is already governed
by the cup product of Schubert classes whose Young diagram is contained
in a small subsquare. This is the content of Theorem \ref{t1}.

We introduce a bilinear pairing on the quantum cohomology ring, which
is just the sum over all coefficients in the Schubert basis expansion
of a product of two classes. The ideal defining the quotient ring
above is just the zero space of the bilinear pairing. In this sense,
the quotient ring describes quantum Schubert calculus modulo {}``numerical
equivalence'' . This results in a simple combinatorial description
of this pairing in the context of symmetric functions in Corollary
\ref{plr}. 

Furthermore, we observe another symmetry related to a parameterization
of points of the quantum cohomology spectrum by Young diagrams as
in \cite{Rie}. The double specialization formula of Stembridge, see
\cite{Ste}, induces a strong symmetry of points and Schubert classes,
which identifies two sets of structures on the quantum cohomology. 

We wish to thank Alexander Postnikov for helpful remarks and conversations.

\section{Complex conjugation and \protect\( \Bbb {Z}/n\protect \) action}

\begin{notation}
Consider \( G(k,n) \), the Grassmannian of \( k- \)planes in \( \Bbb {C}^{n} \)
with tautological bundle \( S \), and the quantum cohomology ring
\( R \) with deformation parameter \( q \) set to \( 1 \). This
ring has the following presentation, see \cite{st}: \[
R=\Bbb {C}[\sigma _{1},\ldots ,\sigma _{k}]/(Y_{n-k+1},\ldots ,Y_{n}+(-1)^{k})\]
\[
l:=n-k,\, \sigma _{i}:=c_{i}(S),\, Y_{i}:=s_{i}(S)\]

We will denote by \( S_{\lambda } \) the Schubert class associated
to the Young diagram \( \lambda =(\lambda _{1},\ldots ,\lambda _{l}) \).
Denote by \( G \) the set of Young diagrams contained in a \( (l,k)- \)rectangle.
The ring R has a vector space basis consisting of Schubert classes
represented by Young diagrams in \( G \). The generators \( \sigma _{i} \)
of \( R \) are represented by diagrams with \( i \) boxes in the
first row and zero boxes in all others. Three point GW-invariants
of cohomology classes will be denoted by \( \left\langle A,B,C\right\rangle  \)
for general classes \( A,B,C \). Poincaré duality on \( R \) will
be denoted by the \( \, \widehat{\hfill }\,  \) -symbol. We use {}``\( * \)''
for the quantum product, and {}``\( \cup  \)'' for the classical
cup product. For standard facts related to the quantum cohomology
ring of Grassmannians we refer to \cite{BCF}. 
\end{notation}
On \( R \), there is a \( \Bbb {Z}/n \) action, given by multiplication
with powers of \( \sigma _{k} \), see \cite[Section 7]{AW}. This
action induces equalities \[
\left\langle \sigma _{k}^{a}*A,\sigma _{k}^{b}*B,\sigma _{k}^{n-a-b}C\right\rangle =\left\langle A,B,C\right\rangle \]
for \( A,B,C\in R, \) see \cite[Proposition 7.2]{AW}. There is an
involution on R, which is just complex conjugation, where we think
of a cohomology class as a complex function on \( \textrm{Spec }R \),
see \cite{Hen}, or \cite{Pos,Rie}, in a different context. This
involution operates on Schubert classes by dualizing the subdiagrams
right and below the Durfee square: \begin{figure}[htbp] \begin{center}

\input{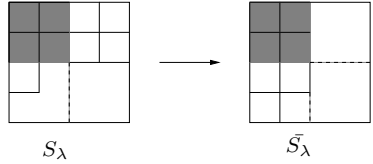}

\caption{Complex conjugation} \label{figure1} \end{center} \end{figure}

For \( \lambda \sim (\lambda _{1},\ldots ,\lambda _{l}), \) writing
\( d_{\lambda } \) for the length of the Durfee square, \( \overline{S_{\lambda }} \)
will be represented by \( (\mu _{1},\ldots ,\mu _{l}) \) with \begin{equation}
\label{eq5}
\mu _{i}=\left\{ \begin{array}{ccc}
d_{\lambda }+k-\lambda _{d_{\lambda }-i+1} & : & i\leq d_{\lambda }\\
d_{\lambda }-\lambda _{l-i+d_{\lambda }+1} & : & i>d_{\lambda }.
\end{array}\right. 
\end{equation}

Algebraically, this involution can be written as

\begin{equation}
\label{defcom}
\overline{S_{\lambda }}=\widehat{S_{\lambda }*C_{P}},
\end{equation}
 where \( C_{P} \) is the class of a point on \( G(k,n). \) 

Complex conjugation yields an important equation of GW-invariants:
\begin{equation}
\label{eq7}
\left\langle A,B,C\right\rangle =\left\langle \widehat{A},\widehat{B},\overline{C}\right\rangle \, \, \, \, A,B,C\in R.
\end{equation}

This equation implies that complex conjugation is a partial inversion:
\( A*\overline{B}=\sum c_{i}C_{i} \), where \( c_{i} \) is the coefficient
of \( A \) in \( C_{i}*B. \) By \ref{eq7}, for any class \( A\in R \),
multiplication by \( A*\overline{A} \) induces a semipositive definite
symmetric endomorphism of the vector space \( R \). Denote by \( P_{1},\ldots ,P_{N} \)
the points of \( \textrm{Spec }R \), then \( A \) is represented
by the vector \( \left( \begin{array}{c}
A(P_{1})\\
\vdots \\
A(P_{N})
\end{array}\right) , \) and \( A*\overline{A} \) will be \( \left( \begin{array}{c}
\left| A(P_{1})\right| ^{2}\\
\vdots \\
\left| A(P_{N})\right| ^{2}
\end{array}\right) . \) 

Because of \[
\overline{\sigma _{k}^{r}*A}=\sigma _{k}^{n-r}*\overline{A},\]
 the map \( A\mapsto A*\overline{A} \) is invariant on the orbits
of \( \Bbb {Z}/n \) action and of course complex conjugation. The
absolute value of Schubert classes should therefore encode its invariant
aspects. 

In the combinatorial setting, we have following property of the map\begin{equation}
\label{mapphi}
\varphi :\, A\mapsto A*\overline{A}:
\end{equation}
 we define for any \( A=\sum c_{\lambda }S_{\lambda }\in R \) \[
\left| A\right| ^{2}=\sum c_{\lambda }\overline{c_{\lambda }},\]
e.g. the square of the hermitian norm with respect to the Schubert
basis. \( \left| A\right| ^{2} \) is the coefficient of 1 in the
expansion of the product \( A*\overline{A} \), since \[
\left\langle A,\overline{A},C_{P}\right\rangle =\sum _{\lambda ,\mu }c_{\lambda }\overline{c_{\mu }}\left\langle S_{\lambda },\overline{S_{\mu }},C_{P}\right\rangle =\sum _{\lambda }c_{\lambda }\overline{c_{\lambda }}\left\langle S_{\lambda },\overline{S_{\lambda }},C_{P}\right\rangle =\sum _{\lambda ,\mu }c_{\lambda }\overline{c_{\mu }},\]
 where we used \ref{defcom}. Since the quantum product is homogenous
mod \( n, \) the quantum ring carries a natural {}``mod \( n \)''
grading. Define \( R_{0}\subseteq R \) to be the subring of homogenous
elements of degree \( 0\textrm{ mod }n, \) and define \( R_{inv} \)
to be the subring of classes invariant under complex conjugation,
that is, the real classes. Then \( \varphi  \) maps Schubert classes
into the subring \( R_{0}\cap R_{inv}. \) For Schubert classes \( S_{\lambda },S_{\mu } \)
the value \[
\left| S_{\lambda }*S_{\mu }\right| ^{2}=S_{\lambda }*S_{\mu }*\overline{S_{\lambda }*S_{\mu }}\mid _{1}\]
 depends only on the values of \( \varphi  \) on \( S_{\lambda },S_{\mu }: \)
The last expression is just the scalar product of the coefficient
vectors of \( \varphi (S_{\lambda }) \) and \( \varphi (S_{\mu }). \)
Below, we will see that a related map controls the sum \( \sum c_{\lambda } \)
in a similar way.

\section{Symmetries of points and Schubert classes}

There is  the following description of the points of \( \textrm{Spec }R \),
see \cite{Rie}. For any \( \lambda \in G \) write \( \lambda =(\lambda _{1},\ldots \, \lambda _{l}) \).
Consider the set \( I \) of all tuples \( I_{\lambda }=(i_{1},\ldots ,i_{l}) \)
indexed by \( \lambda \in G, \) where \( i_{j}=\frac{l+1}{2}+\lambda _{j}-j. \)
We fix \( \xi  \) to be any primitive \( n \)-th root of unity,
and, in case \( l \) even, \( \xi ^{\frac{1}{2}} \) to be any primitive
\( 2n \) -th root of unity. Then to any \( I_{\lambda } \) we can
associate a point \( P_{\lambda }\in \textrm{Spec }R\hookrightarrow \Bbb {C}^{k} \)
by the following map:\[
I_{\lambda }=(i_{1},\ldots ,i_{l})\mapsto (h_{1}(\xi ^{i_{1}},\ldots ,\xi ^{i_{l}}),\ldots ,h_{k}(\xi ^{i_{1}},\ldots ,\xi ^{i_{l}})),\]
 where the \( h_{i} \) are the full homogenous symmetric polynomials
of degree \( i \) in \( l \) variables. Abusing notation, we will
sometimes write \( P_{S} \) for the point corresponding to any Schubert
class \( S, \) so \( P_{1} \) is the point associated to the empty
diagram. \( S_{\lambda }(P_{1}) \) is invariant under complex conjugation,
therefore real, and also positive, see \cite[Theorem 8.4]{Rie}. 

\begin{prop}
\label{p2}(Double specialization formula) Consider the matrix \( M:=(S_{\lambda }(P_{\mu }))_{\lambda ,\mu } \).
M is symmetric modulo division by the \( P_{1} \), that is \begin{equation}
\label{id1}
\frac{S_{\lambda }(P_{\mu })}{S_{\lambda }(P_{1})}=\frac{S_{\mu }(P_{\lambda })}{S_{\mu }(P_{1})}.
\end{equation}
 
\end{prop}
\begin{proof}
This identity is an easy consequence of the definition of Schur polynomials,
see \cite[Exercise 7.32a]{Stan}:

As a symmetric function in \( k \) variables \( x_{1},\ldots ,x_{k} \),
the Schur polynomial can be defined as\[
S_{\lambda }=\frac{\textrm{det }\left( x_{i}^{\lambda _{j}+k-j}\right) _{1\leq i,j\leq k}}{\textrm{det }\left( x_{i}^{k-j}\right) _{1\leq i,j\leq k}},\]
 see \cite{mac}, I (3.1). In this sense, \( S_{\lambda }(P_{\mu }) \)
will be \( S_{\lambda }(\xi ^{I_{\mu }}). \) Proving the identity
\ref{id1}, because of cancellation, it is enough to calculate the
left hand side with \( I'_{\mu }=(\mu _{1}+k-1,\mu _{2}+k-2,\ldots ,\mu _{k}) \).
By \cite[Lemma 4.4]{Rie}, we can assume \( l\leq k \). Filling up
the Young diagrams with zero rows we can assume \( l=k \). We need
to prove:\[
\frac{S_{\lambda }(\xi ^{I'_{\mu }})}{S_{\lambda }(\xi ^{I_{1}'})}=\frac{S_{\mu }(\xi ^{I_{\lambda }'})}{S_{\mu }(\xi ^{I_{1}'})}.\]
\[
\frac{S_{\lambda }(\xi ^{I'_{\mu }})}{S_{\lambda }(\xi ^{I_{1}'})}=\frac{\frac{\textrm{det }\left( \xi ^{\mu _{i}+k-i+\lambda _{j}+k-j}\right) _{i,j}}{\textrm{det }\left( \xi ^{\mu _{i}+k-i+k-j}\right) _{i,j}}}{\frac{\textrm{det }\left( \xi ^{k-i+\lambda _{j}+k-j}\right) _{i,j}}{\textrm{det }\left( \xi ^{k-j+k-j}\right) _{i,j}}}=\frac{\textrm{det }\left( \xi ^{\mu _{i}+k-i+\lambda _{j}+k-j}\right) \cdot \textrm{det }\left( \xi ^{k-i+k-j}\right) }{\textrm{det }\left( \xi ^{k-i+\lambda _{j}+k-j}\right) \cdot \textrm{det }\left( \xi ^{\mu _{i}+k-i+k-j}\right) }.\]
 We observe that the expression on the right is symmetric in \( \lambda  \)
and \( \mu  \). 
\end{proof}
\begin{rem}
The identity \ref{id1} is a special case of a double specialization
formula of Stembridge, see \cite[Therorem 7.4]{Ste}. This paper also
includes reference to the Macdonald polynomial conjecture, which is
a further generalization. 
\end{rem}
\begin{example}
\( n=4,\, \, k=2,\, \, \xi =i,\, \xi ^{\frac{1}{2}}=\frac{i+1}{\sqrt{2}} \)
:\[
\begin{array}{ccccccc}
 & S_{0,0} & S_{1,0} & S_{2,0} & S_{1,1} & S_{2,1} & S_{2,2}\\
P_{0,0} & 1 & \sqrt{2} & 1 & 1 & \sqrt{2} & 1\\
P_{1,0} & 1 & 0 & i & -i & 0 & -1\\
P_{2,0} & 1 & \sqrt{2}i & -1 & -1 & -\sqrt{2}i & 1\\
P_{1,1} & 1 & -\sqrt{2}i & -1 & -1 & \sqrt{2}i & 1\\
P_{2,1} & 1 & 0 & -i & i & 0 & -1\\
P_{2,2} & 1 & -\sqrt{2} & 1 & 1 & -\sqrt{2} & 1
\end{array}\]

After division of the row corresponding to \( P_{0,0} \) we get a
symmetric matrix: \[
\begin{array}{ccccccc}
 & S_{0,0} & S_{1,0} & S_{2,0} & S_{1,1} & S_{2,1} & S_{2,2}\\
P_{0,0} & 1 & 1 & 1 & 1 & 1 & 1\\
P_{1,0} & 1 & 0 & i & -i & 0 & -1\\
P_{2,0} & 1 & i & -1 & -1 & -i & 1\\
P_{1,1} & 1 & -i & -1 & -1 & i & 1\\
P_{2,1} & 1 & 0 & -i & i & 0 & -1\\
P_{2,2} & 1 & -1 & 1 & 1 & -1 & 1
\end{array}\]

\end{example}
We have the following result for the values of the \( S_{\lambda } \)
on \( P_{\sigma _{k}} \), or the value of \( \sigma _{k} \) on points
\( P_{\lambda } \):

\begin{lem}
\label{lem1}\( \frac{S_{\lambda }(P_{\sigma _{k}})}{S_{\lambda }(P_{1})}=\sigma _{k}(P_{\lambda })=\xi ^{d}, \)
where \( d \) is the degree of \( S_{\lambda } \) modulo \( n \). 
\end{lem}
\begin{proof}
We note that, by definition, for any \( I_{\lambda }=(i_{1},\ldots ,i_{l}) \)
the sum \( \sum i_{j} \) is just the degree of \( S_{\lambda } \).
Now, taking \( e_{i} \) to be the \( i \)-th elementary symmetric
polynomial and \( \lambda ^{t}=(\lambda ^{t}_{1},\ldots ,\lambda _{k}^{t}) \)
the transposed Young diagram, \( h_{k}(\xi ^{I_{\lambda }})=e_{k}(\xi ^{I_{\lambda ^{t}}}), \)
see \cite[Proof of Lemma 4.4]{Rie}, and \( e_{k}(\xi ^{I_{\lambda ^{t}}})=\xi ^{\lambda ^{t}_{1}}\cdots \xi ^{\lambda _{k}^{t}}=\xi ^{\left| I_{\lambda ^{t}}\right| }=\xi ^{\left| I_{\lambda }\right| }. \) 
\end{proof}
One immediate consequence of this Lemma is that the number of orbits
under the \( \Bbb {Z}/n \) action is the same as the number of Schubert
classes in \( R_{0} \), the subring of classes of degree zero modulo
\( n \):

\begin{prop}
\label{orbits}\( \textrm{dim}_{\Bbb {C}}\, R_{0}=\#\{\textrm{Orbits under the }\Bbb {Z}/n\, \textrm{action}\} \)
\end{prop}
\begin{proof}
The number of orbits under the \( \Bbb {Z}/n \) action is equal to
the dimension of \( R/(\sigma _{k}-1) \), since the relation \( (\sigma _{k}-1) \)
identifies exactly Schubert classes in the same orbit. Since we are
dividing out a hyperplane, \( \textrm{dim }R/(\sigma _{k}-1) \) will
be the number of points \( P_{\lambda }\in \textrm{Spec }\, R \)
with \( \sigma _{k}(P_{\lambda })=1. \) The statement now follows
from the previous lemma.
\end{proof}
The equivalent of the \( \Bbb {Z}/n \) action on points is multiplication
with roots of unity. Defining \( \Xi =(\xi ,\xi ^{2},\ldots ,\xi ^{k}) \),
the action is induced by component wise multiplication of the vector
\( P_{\lambda }\sim (\sigma _{1}(P_{\lambda }),\ldots ,\sigma _{k}(P_{\lambda })) \)
with powers of \( \Xi  \).

\section{Quantum cohomology modulo symmetries}

\begin{figure}[!h] \begin{center}

\input{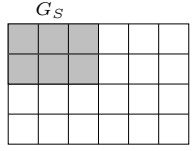}

\caption{The \( l\times k \) box with subrectangle \( G_{S} \) } \label{figure3} \end{center} \end{figure}

\begin{defn}
(see Figure \ref{figure3}) For the cohomology ring \( G \) of the
Grassmannian \( G(k,n) \), with \( l:=n-k, \) define the subvectorspace
\( G_{S} \) to be spanned by the classes associated to Young diagrams
contained in the \( [\frac{k}{2}]\times [\frac{l}{2}] \) rectangle,
where \( [\, ] \) means rounding down. Abusing notation, we will
sometimes consider \( G_{S} \) to be the set of Young diagrams contained
in this rectangle.

We define a map from \( G_{S}\times G_{S}\mapsto G \) by assigning
two Schubert classes \( S_{\lambda } \) and \( S_{\mu } \) the class
\( S_{\lambda }(S_{\mu }) \) defined by embedding the Young diagram
of \( S_{\mu } \) into the diagram of \( \overline{S_{\lambda }} \)
in the following way:\[
n=8,\, k=4,\, \lambda =\yng (2,1),\, \mu =\yng (1,1),\, \overline{\lambda }=\yng (3,1,1)\]
\[
S_{\lambda }(S_{\mu })\cong \yng (3,2,2)\]
 Such an embedding is always possible, since for \( S_{\lambda }\in G_{S}, \)
\( \overline{S_{\lambda }} \) always has an empty \( [\frac{k}{2}]\times [\frac{l}{2}] \)
rectangle at the outer corner of the Durfee square. Numerically, with
\( d_{\lambda } \) the length of the Durfee square, and \( \lambda \cong (\lambda _{1},\ldots  \),\( \lambda _{l}),\, \mu \cong (\mu _{1},\ldots ,\mu _{l}) \),
\( S_{\lambda }(S_{\mu }) \) is given by \( (\nu _{1},\ldots ,\nu _{l}) \)
\begin{equation}
\label{eq1}
\nu _{i}=\left\{ \begin{array}{ccc}
d_{\lambda }+k-\lambda _{d_{\lambda }-i+1} & : & i\leq d_{\lambda }\\
\mu _{i-d_{\lambda }}+d_{\lambda }-\lambda _{l-i+d_{\lambda }+1} & : & i>d_{\lambda }.
\end{array}\right. 
\end{equation}

\end{defn}

\begin{rem}
\label{rem1}The subsquare \( G_{S} \) is the biggest subrectangle
with the property that quantum product is just ordinary cup product
without any {}``cutting off''.
\end{rem}
\begin{prop}
\label{p3}All Schubert classes of degree zero mod n invariant under
complex conjugation are of type \( T_{\alpha }(T_{\alpha }) \) for
some \( \alpha \in G_{S}. \)
\end{prop}
\begin{proof}
Take a Schubert class \( S_{\lambda }\in G_{S} \). We claim that
the Young diagram \( \overline{S_{\lambda }}(S_{\lambda }) \) is
invariant. This follows from a direct calculation using \ref{eq1}
and \ref{eq5}.\\
Now let \( S_{\lambda } \) be a degree \( 0 \) invariant Schubert
class. By \ref{eq5}, \[
\textrm{deg }\overline{S_{\lambda }}=nd_{\lambda }-\textrm{deg }S_{\lambda },\]
 so for \( S_{\lambda } \) we have \( nd_{\lambda }=2\, \textrm{deg S}_{\lambda }=2\, na \)
for some \( a\in \Bbb {N}. \) The length of the Durfee square is
even, so we can reverse the construction. With \( S_{\lambda }=(\lambda _{1},\ldots ,\lambda _{l}), \)
the corresponding class \( S\in G_{S} \) is given by \( \mu =(\mu _{1},\ldots ,\mu _{\left[ \frac{l}{2}\right] },0,\ldots ,0) \),\[
\mu _{i}=\lambda _{d_{\lambda }/2+i}-d_{\lambda }/2,\, \, \, i\leq \left[ \frac{l}{2}\right] .\]
 It follows from \ref{eq5} that \( S_{\lambda } \) is already determined
by the \( \mu _{i} \), and that \( \mu  \) is contained in \( G_{S} \):
the subdiagrams right and below of the Durfee square of \( \mu  \)
correspond to lower halves of selfdual Young diagrams, so they can
not extend over more than half of the length or height of the rectangle.
Also, the Durfee square itself has length \( \frac{d_{\lambda }}{2}. \) 
\end{proof}
We prove the following important description for the product of two
classes in \( G_{S} \):

\begin{prop}
\label{p1}For \( S_{\lambda },S_{\mu }\in G_{S}, \)\begin{equation}
\label{f3}
\left\langle \overline{S_{\lambda }},S_{\mu }\widehat{,T_{\alpha }(T_{\beta })}\right\rangle =\sum _{T_{\nu }\in G_{S}}\left\langle \widehat{S_{\lambda }},T_{\alpha },T_{\nu }\right\rangle \left\langle T_{\nu },T_{\beta },\widehat{S_{\mu }}\right\rangle ,
\end{equation}
 and \( \overline{S_{\lambda }}*S_{\mu }=\sum _{\alpha ,\beta }\left\langle \overline{S_{\lambda }},S_{\mu }\widehat{,T_{\alpha }(T_{\beta })}\right\rangle T_{\alpha }(T_{\beta }). \) 
\end{prop}
\begin{proof}
We will use the description of quantum product by rim hook reduction
as in \cite{BCF}, see also \cite{Buc} for a simple proof. All diagrams
with nontrivial reduction will have a structure as indicated in Figure
\ref{figure2}, I.

\begin{figure}[htbp] \begin{center}

\input{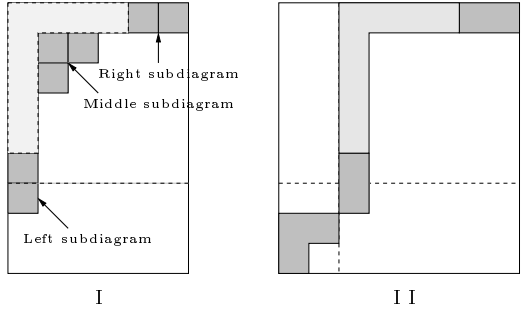}\label{fig2}

\caption{Skew diagrams associated to \( \overline{S_{\lambda }}\cup S_{\mu } \) and \( \overline{S_{\lambda }}*S_{\lambda } \)} \label{figure2} \end{center} \end{figure}

We will now enumerate all semistandard skew tableaux appearing in
\( \overline{S_{\lambda }}\cup S_{\mu } \) whose associate word is
a  reverse lattice word and whose reduction is of shape \( T_{\alpha }(T_{\beta }) \)
for fixed \( \alpha ,\beta  \). Note that such diagrams must have
\( \beta  \) as middle subdiagram, since any removable rim hook must
extend over the width and length of \( G_{S} \), such that the reduction
will again have a shape similar to Figure \ref{figure2}, with \( \beta  \)
inscribed in the middle. Now take the set of all diagrams reducing
to \( T_{a}(T_{\beta }) \) by the rim hook reduction. For any such
fixed diagram, the number of reverse lattice word completions is the
same as for the skew diagram obtained by exchanging the position of
the middle and left subdiagram as in Figure \ref{figure2}, II , by
reason of commutativity and associativity of the product in the tableaux
ring, see \cite[§5.1]{Ful}. It is easy to enumerate theses skew tableaux:
looking at the set of skew diagrams with content \( T_{\nu } \) in
the right and middle subdiagram we find that this number is just \( \left\langle \overline{S_{\lambda }},T_{\nu },\widehat{\overline{T_{\alpha }}}\right\rangle =\left\langle S_{\lambda },\widehat{T_{\nu }},\overline{T_{\alpha }}\right\rangle =\left\langle \widehat{S_{\lambda }},T_{\nu },T_{\alpha }\right\rangle , \)
where we used \ref{eq7}, since the rim hook reduction is independent
of the subdiagram inscribed in the middle. Now, for a fixed completion
of the middle and right subdiagram of content \( T_{\nu } \), we
can enumerate all possible completions of the left one. By one of
the definitions of the Littlewood-Richardson rule, see \cite[§5.2, Corollary 2]{Ful},
we find that this number is indeed \( \left\langle T_{\nu },T_{\beta },\widehat{S_{\mu }}\right\rangle  \).
Summing up, we get \ref{f3}. The argument above also shows that rim
hook reduction will only result in diagrams of shape \( T_{\alpha }(T_{\beta }) \)
for some \( \alpha ,\beta \in G_{S} \). 
\end{proof}
\begin{rem}
Using \ref{eq7}, the formula \ref{f3} can be interpreted as a restricted
associativity equation for the product \( S_{\lambda }*\overline{T_{\alpha }}*T_{\beta } \).
In this sense, the class \( T_{\alpha }(T_{\beta }) \), which appears
in \( \overline{T_{\alpha }}*T_{\beta } \), corresponds to restricting
the associativity equation to classes \( T_{\nu } \) in \( G_{S}. \) 
\end{rem}
\begin{example}
Let \( n=12,\, k=6,\, \lambda =\tiny {\yng (3,2,1)},\, \mu =\tiny {\yng (2,2,1)}, \)
\( \alpha =\tiny {\yng (2,2)},\, \beta =\tiny {\yng (2,1)} \). The
coefficient of \( T_{\alpha }(T_{\beta })\cong \tiny {\yng (6,6,4,3,2,2)} \)
in \( \overline{S_{\lambda }}*S_{\mu } \) is equal to the number
of pairs of semistandard skew tableaux on the shapes \( \lambda /\alpha  \)
and \( \mu /\beta  \) of same content, satisfying the reverse lattice
word condition, as in the figure below:\[
\young (\hfil \hfil 1,\hfil \hfil ,2)\, \, \, \, \, \, \, \, \, \, \, \, \, \, \, \, \, \young (\hfil \hfil ,\hfil 1,2).\]
In our example, \( \left\langle \overline{S_{\lambda }},S_{\mu }\widehat{,T_{\alpha }(T_{\beta })}\right\rangle  \)
will be \( 2. \) 

Associate to \( S_{\lambda }\in G_{S} \) the symmetric matrix \[
A_{\lambda }:=\left( \left\langle \widehat{S_{\lambda }},T_{\alpha },T_{\beta }\right\rangle \right) _{\alpha ,\beta \in G_{S}}.\]
 The product \( \overline{S_{\lambda }}*S_{\mu } \) is determined
by the matrix product \( A_{\lambda }*A_{\mu }^{t}. \)

The strong combinatorial description of the product \( \overline{S_{\lambda }}*S_{\mu } \)
should be helpful for further understanding the {}``unlimited''
cup product \( S_{\lambda }*S_{\mu }=S_{\lambda }\cup S_{\mu }. \)
The relationship between \( \overline{S_{\lambda }}*S_{\mu } \) and
\( S_{\lambda }*S_{\mu } \) reflects a symmetry between adding and
removing skew diagrams in a quantum sense, see \ref{eq7}. In the
following we will see that \( S_{\lambda } \) and \( \overline{S_{\lambda }} \)
display an identical {}``numerical'' behaviour. 
\end{example}
\begin{defn}
Define the class \( P:=\sum _{\lambda }S_{\lambda } \). The ideal
Ann \( P \) should be considered as consisting of relations \( (S_{\lambda }-\phi (S_{\lambda })) \)
where \( \phi  \) is a {}``symmetry map''. We will also consider
the bilinear pairing \( \left\langle A,B\right\rangle _{s}=\sum _{T_{\lambda }\in G}<A,B,T_{\lambda }>=\left\langle A,B,P\right\rangle , \)
which is just the sum over all coefficients appearing in the Schubert
basis representation of products of classes \( A,B\in R. \) \( \textrm{Ann }P \)
is obviously the zero space of \( \left\langle \cdot ,\cdot \right\rangle _{s} \)
and can be considered to be the ideal of relations of numerically
equivalent classes with respect to \( \left\langle \cdot ,\cdot \right\rangle _{s}. \) 
\end{defn}
Taking \( \sigma _{0}:=1 \), we wish to show the following 

\begin{thm}
\label{t1}\( \textrm{Ann }P\cong (\{\sigma _{k-i}-\sigma _{i}\}_{i}) \),
and \( R_{P}:=R/\textrm{Ann }P \) as a vector space is isomorphic
to \( G_{S}. \) 
\end{thm}
\begin{rem}
\label{abbphi}We have \( (\{\sigma _{k-i}-\sigma _{i}\})=(\{\overline{\sigma _{i}}-\sigma _{i}\},\sigma _{k}-1)\supsetneq \{S_{\lambda }-\overline{S_{\lambda }},S_{\lambda }-\sigma ^{i}_{k}S_{\lambda }\, :\, S_{\lambda }\in R\} \),
which shows that \( \textrm{Ann }P \) identifies all the orbits under
the \( \Bbb {Z}/n \) action and complex conjugation. 

Let \( V \) be a vector space with basis \( \{v_{\lambda }\} \)
indexed by the Schubert classes in \( G_{S} \). Then the map\begin{equation}
\label{mappsi}
\begin{array}{ccc}
\psi :\, G_{S} & \rightarrow  & V\\
S_{\lambda } & \mapsto  & \sum _{\nu \in G_{S}}\sum _{\lambda \in G_{S}}\left\langle \widehat{S},T_{\lambda },T_{\mu }\right\rangle v_{\nu }
\end{array}
\end{equation}
is an isomorphism of vector spaces, since the corresponding matrix
can be arranged to be triangular with \( 1 \) on the diagonal: the
coefficient of \( v_{\nu } \) in \( \psi (S_{\lambda }) \) will
only be nonzero if the diagram \( \nu  \) is contained in \( \lambda  \),
and \( 1 \) if \( \nu =\lambda . \) The bilinear form \( \left\langle -,-\right\rangle _{s} \)
simply becomes the standard scalar product under this map. Since,
by the Theorem, \( \left\langle -,-\right\rangle _{s} \) descends
to \( R_{P} \), we deduce the following Corollary, see Remark \ref{rem1}:
\end{rem}
\begin{cor}
\label{plr}Consider the Littlewood-Richardson expansion of the product
of two Schur polynomials in the ring of symmetric functions: \[
S_{\lambda }*S_{\mu }=\sum c^{v}_{\lambda ,\mu }S_{\nu }.\]
 The sum \( \sum _{\nu }c_{\lambda ,\mu }^{\nu } \)is identical to
\[
\sum _{\alpha ,\beta ,\nu }c_{\alpha ,\nu }^{\lambda }\cdot c_{\beta ,\nu }^{\mu }.\]

\end{cor}
\begin{proof}
(of the Theorem) Define \( I:=(\{\sigma _{k-i}-\sigma _{i}\}_{i}). \)We
first show that \( I\subseteq \textrm{Ann }P \). It is obvious that
\( \sigma _{k}-1\in \textrm{Ann P} \). We will show \( \overline{\sigma _{i}}-\sigma _{i}\in \textrm{Ann }P, \)
which by \ref{eq2} implies \( \sigma _{k-i}-\sigma _{i}\in \textrm{Ann }P. \)
By an application of complex conjugation we find \( \sigma _{i}*P\mid \widehat{_{S_{\lambda }}}=\overline{\sigma _{i}}*P\mid \widehat{_{S_{\lambda }}} \)
: \[
\begin{array}{ccl}
\left\langle \sigma _{i},P,S_{\lambda }\right\rangle  & = & \sum \left\langle \sigma _{i},S_{\lambda },T_{i}\right\rangle =\sum \left\langle \sigma _{i},S_{\lambda },T_{i}\right\rangle ^{2}=\sigma _{i}*S_{\lambda }*\overline{S_{\lambda }}\cup \widehat{\sigma _{i}}\\
 & = & \sum \left\langle \sigma _{i},\overline{\sigma _{i}},T_{i}\right\rangle \left\langle S_{\lambda },\overline{S_{\lambda }},T_{i}\right\rangle =\sum \left\langle \overline{\sigma _{i}},S_{\lambda },T_{i}\right\rangle ^{2}\\
 & = & \sum \left\langle \overline{\sigma _{i}},S_{\lambda },T_{i}\right\rangle =\left\langle \overline{\sigma _{i}},P,S_{\lambda }\right\rangle ,
\end{array}\]

using \ref{plr}, \ref{eq2} and associativity of quantum multiplication
together with the fact that \( \left\langle \sigma _{i},S_{\lambda },T_{i}\right\rangle  \)
is either \( 1 \) or \( 0 \). We observe that the classes in \( G_{S} \)
have no relations in \( \textrm{Ann }P \): this can be translated
into the fact that \( \left\langle \cdot ,\cdot \right\rangle _{s}=\left\langle \cdot ,\cdot ,P\right\rangle  \)
is nondegenerate on \( G_{S}, \) which follows from the discussion
after the Theorem. 

Now it is enough to show that the dimension of \( R/I \) does not
exceed the number of Schubert classes in \( G_{S}. \) We observe
that \( \textrm{Spec }R/I \) consists of all real points \( P_{\lambda }\in \textrm{Spec }R \)
with \( \sigma _{k}(P_{\lambda })=1 \), since \begin{equation}
\label{eq2}
\overline{\sigma _{i}}=\sigma _{k-i}*\sigma _{k},
\end{equation}
 so the relations \( \sigma _{k-i}-\sigma _{i} \) guarantee that
the imaginary part of any \( P_{\lambda } \) with \( \sigma _{k}(P_{\lambda })=1 \)
vanishes. But these points, by Proposition \ref{p2} and Lemma \ref{lem1},
correspond 1-1 to Schubert classes of degree zero mod n which are
invariant under complex conjugation. By Proposition \ref{p3}, these
are parameterized by Schubert classes in \( G_{S}. \) We have proved\[
\#\{P_{\lambda }\in \textrm{Spec}\, R/I\}\leq \#\{S_{\lambda }\in G_{S}\},\]
 which implies the result. 
\end{proof}
\begin{rem}
1. The last statement in the previous proof can be seen directly by
proving that every monomial in \( R \) can be represented mod \( I \)
by monomials coming from classes in \( G_{S}. \) This is done by
carrying out the identification of both generators \( \sigma _{i} \)
\( \sim \sigma _{k-i} \) and dual generators \( \sigma ^{i}\sim \sigma ^{l-i}. \) 

2. The points of \( \textrm{Spec }R_{P} \) are exactly the \( P_{\lambda } \)
for which \( S_{\lambda } \) is invariant of degree zero mod n. \( R_{P} \)
is the coordinate ring of these points and the quotient map \( R\mapsto R_{P} \)
is the restriction map. All elements of \( R_{P} \) are real functions. 

3. It is obvious that all the Schubert classes in the orbit of one
\( S_{\lambda }\in G_{S} \) are mapped to this class under the composition
\[
R\rightarrow R_{P}\cong G_{S}.\]
 In general, a Schubert class will be mapped to a nonpositive linear
combination of Schubert classes in \( G_{S}, \) as in the following
example:\\
 \( n=8,\, k=4\, \Rightarrow \, \, G_{S}\sim G(2,4) \).\\
Take \[
\lambda =\yng (3,2,1).\]
 In \( G_{S} \), \( \left[ \lambda \right]  \) is represented by
\[
\yng (2,2)+\yng (2,0)+\yng (1,1)-2.\]

Under the isomorphism \( \psi  \) of \ref{mappsi}, \( \left[ \lambda \right]  \)
is represented by \[
\yng (2,2)+\yng (2,1)+2\yng (2,0)+2\yng (1,1)+3\yng (1,0)+1.\]
 Note that this is a positive linear combination.  
\end{rem}

\section{Summary and Outlook}

In this paper we observed the following set of dualities:\[
\begin{array}{rcl}
\{\textrm{Orbits under }\Bbb {Z}/n\textrm{ action}\} & \approx  & R_{0}\\
\{\textrm{Orbits under }\Bbb {Z}/n\textrm{ action and comp}.\textrm{ conj}.\} & \approx  & R_{0}\cap R_{inv}\\
R_{P} & \approx  & \{S_{\lambda }\in R_{0}\cap R_{inv}\}.
\end{array}\]
 The symmetry of Proposition \ref{p2} hints at a stronger relation
beyond the equality of dimensions as vector spaces. This is especially
obvious in the third case, where, as in the case of \( R, \) we find
a 1-1 correspondence between points of \( \textrm{Spec }R_{P} \)
and a {}``Schubert basis'' of \( R_{P}. \) The ring \( R_{P} \)
controls quantum Schubert calculus modulo numerical equivalence, which
implies that the orbits under the actions behave numerically identical.
We expect that the composition \( R\rightarrow R_{P}\stackrel{\psi }{\rightarrow }V \)
results in vectors with positive components for all Schubert classes.
A positive combinatorial formula for this map would suggest a new
rule for constructing the quantum product. This is especially interesting
with regards to the positivity of the GW invariants, which is evident
geometrically but not as clear algebraically. 

It would be rewarding to find an explicit relationship between products
\( S_{\lambda }*S_{\mu } \) and \( \overline{S_{\lambda }}*S_{\mu } \),
since this would imply a better understanding of the classical cup
product, e.g. with regards to the question as to what classes appear
with nonzero coefficient. Also, there might be a natural generalization
to the quantum product. 

A part of the results for the Grassmannian should be extendable to
other homogenous varieties. For example, the variety of full flags
in \( \Bbb {C}^{n} \) carries a \( \Bbb {Z}/n \) action, see \cite{Pos2},
which acts by shifting the word of a permutation. Computations indicate
that the zero space of \( \left\langle \cdot ,\cdot \right\rangle _{s} \)
is generated by this action. In analogy to the Grassmannian, the quotient
ring would have a vector space basis consisting e.g. of Schubert classes
whose permutations fix \( n. \) The quantum product should be nondegenerate
on this subset, that is, the linear map \( \left\langle S_{\omega },\cdot \right\rangle _{s} \)
determines the Schubert class. In analogy to Proposition \ref{orbits},
the numer of orbits, \( (n-1)!, \) is identical to the number of
Schubert classes of degree zero mod \( n, \) because the number of
permutations with \( t\textrm{ mod }n \) inversions is indeed \( (n-1)! \),
see for example \cite[5.1.1, Exercise 13]{kn}.

\end{document}